\documentclass[graybox]{svmult}

\usepackage{mathptmx}       
\usepackage{helvet}         
\usepackage{courier}        
\usepackage{type1cm}        
\usepackage{makeidx}         
\usepackage{graphicx}        
\usepackage{multicol}        
\usepackage[bottom]{footmisc}

\usepackage{amsmath,amsfonts,amssymb}
\usepackage{booktabs}
\usepackage{epstopdf}
\usepackage{textcomp}

\makeindex             

\begin{document}

\title*{Space-Time NURBS-Enhanced Finite Elements for Solving the Compressible Navier-Stokes Equations}
\titlerunning{Space-Time NEFEM for Solving the Compressible Navier-Stokes Equations} 
\author{Michel Make, Norbert Hosters, Marek Behr and Stefanie Elgeti}
\institute{Michel Make \at Chair for Computational Analysis of Technical Systems, RWTH-Aachen University, Schinkelstr. 2, 52062 Aachen, Germany, \email{make@cats.rwth-aachen.de}
\and Norbert Hosters \at Chair for Computational Analysis of Technical Systems, RWTH-Aachen University, Schinkelstr. 2, 52062 Aachen, Germany,  \email{hosters@cats.rwth-aachen.de}
\and Marek Behr \at Chair for Computational Analysis of Technical Systems, RWTH-Aachen University, Schinkelstr. 2, 52062 Aachen, Germany, \email{behr@cats.rwth-aachen.de}
\and Stefanie Elgeti \at Chair for Computational Analysis of Technical Systems, RWTH-Aachen University, Schinkelstr. 2, 52062 Aachen, Germany,\email{elgeti@cats.rwth-aachen.de}
\\\\
\and \textit{NOTICE: This is a pre-print of an article submitted for publication in Lecture Notes in Computational Science and Engineering. Changes resulting from the publishing process, such as editing, corrections, structural formatting, and other quality control mechanisms may not be reflected in this document. Changes may have been made to this work since it was submitted for publication.
}}
%
\maketitle

\abstract{This article considers the NURBS-Enhanced Finite Element Method (NEFEM) applied to the compressible Navier-Stokes equations. NEFEM, in contrast to conventional finite element formulations, utilizes a NURBS-based computational domain representation. Such representations are typically available from Computer-Aided-Design tools. Within the NEFEM, the NURBS boundary definition is utilized only for elements that are touching the domain boundaries. The remaining interior of the domain is discretized using standard finite elements. Contrary to isogeometric analysis, no volume splines are necessary. 
\newline\indent
The key technical features of NEFEM will be discussed in detail, followed by a set of numerical examples that are used to compare NEFEM against conventional finite element methods with the focus on compressible flow.}

\section{Introduction}
\label{sec:introduction}

Geometries in engineering applications are commonly designed with the use of Computer-Aided-Design (CAD) tools. In general, these tools utilize Non-Uniform Rational B-Splines (NURBS) to accurately represent complex geometric domains by means of surface splines. When performing numerical analysis on such domains, it is common practice to first discretize the domain into finite sub-domains or elements. This discretization process, typically results in loss of the exact geometry. 

An alternative approach, known as Isogeometric Analysis (IGA), was proposed in \cite{hughes2005}. The key idea of IGA is to use the NURBS basis functions not only for the geometric representation, but also for the numerical solution itself. By doing so, numerical analysis can be applied to the CAD model directly without the loss of geometric accuracy caused by discretizing the computational domain. Numerical analysis of fluid flow problems, however, commonly involves complex three-dimensional volume domains. Parametrizing such domains using closed volume splines can be challenging. 

An alternative was proposed in \cite{sevilla2008}, and further extended for space-time finite elements and free-surface flows in \cite{stavrev2016}. This method was then modified for interface-coupled problems in \cite{hosters2018}. This approach suggests to use standard finite elements in the interior of the computational domain supplemented with so-called NURBS-enhanced finite elements along domain boundaries. These elements make use of NURBS to accurately represent complex geometries. The NURBS-Enhanced Finite Element Method (NEFEM) allows for maintaining as much as possible the proven computational efficiency of standard finite element methods, while utilizing the accurate geometric representation provided by the NURBS.

In this work, we apply NEFEM to supersonic flow problems. For this type of problems, accurate geometry representation can be important, especially due to the presence of shock waves and their interaction with solid walls. 

\section{Quasi-Linear Form of the Navier-Stokes equations}
\label{sec:governing-equations}
Before presenting the NEFEM concept, first the governing Navier-Stokes equations are presented. For this, let $\Omega_t\subset \mathbb{R}^{n_{sd}}$ and $t \in (0,T)$ be the spatial and temporal domains respectively, and let $\Gamma_t$ denote the boundary of  $\Omega_t$. Then the model problem, written as a generalized advective-diffusive system, is given by:
\begin{align}
\frac{\partial \mathbf{U}      }{\partial t} + \frac{\partial \mathbf{F}_i    }{\partial x_i} - \frac{\partial \mathbf{E}_i    }{\partial x_i} &= \mathbf{0} \quad \quad\mbox{on $\Omega_t \quad \forall\:t \in (0,T)$},\label{eq:NS-syseq}
\end{align}
where $\mathbf{U} = (\rho, \rho u_1, \rho u_2, \rho u_3, \rho e)^T$ is the solution vector. $\rho$, $u_i$, and $e$ represent the density, velocity components, and total energy per unit mass respectively. For the three-dimensional case, the Euler and viscous flux vectors $\mathbf{F}_i$ and $\mathbf{E}_i$ are defined as: 
\begin{equation}
\mathbf{F}_i = \begin{pmatrix} u_i \rho \\  
                      u_i \rho u_1 + \delta_{i1} p \\
                      u_i \rho u_2 + \delta_{i2} p \\
                      u_i \rho u_3 + \delta_{i3} p \\
                      u_i (\rho e +  p)  \end{pmatrix}, \\
                      \qquad
 \mathbf{E}_i = \begin{pmatrix} 0 \\  
                      \tau_{i1} \\
                      \tau_{i2} \\
                      \tau_{i3} \\
                      -q_i + \tau_{ik} u_k \end{pmatrix},                
\end{equation}
where $q_i$, $\tau_{ik}$, $p$, and $\delta$ represent the heat flux, viscous stress tensor, pressure, and the Kronecker delta respectively. The boundary and initial conditions are given by:
\begin{align}
\mathbf{U} \cdot \mathbf{e}_d &= g_d  \quad \mbox{on $(\Gamma_t)_{g_d}, \quad d = 1...\;n_{dof} \:\forall\:t \in (0,T)$},\\
(n_i\mathbf{E}_i) \cdot \mathbf{e}_d &= h_d \quad\mbox{on $(\Gamma_t)_{h_d}, \quad d = 1...\;n_{dof} \:\forall\:t \in (0,T)$},\\
\mathbf{U}(\mathbf{x},0) &= \mathbf{U}_0 (\mathbf{x})         \quad \mbox{on $\Omega_0$},
\end{align}
where $(\Gamma_t)_{g_d}$ and $(\Gamma_t)_{h_d}$ are the subsets of $\Gamma_t$, $\mathbf{e}_d$ is a basis in $\mathbb{R}^{n_{dof}}$, and $n_{dof}$ is the number of degrees of freedom. 
The quasi-linear form of Equation \eqref{eq:NS-syseq} is written as:
\begin{align}
\frac{\partial \mathbf{U}      }{\partial t} + \mathbf{A}_i\frac{\partial \mathbf{U}     }{\partial x_i} - \frac{\partial                  }{\partial x_i} \left( \mathbf{K}^h_{ij} \frac{\partial \mathbf{U}     }{\partial x_j} \right) = \mathbf{0} \quad \mbox{on $\Omega_t \quad \forall\:t \in (0,T)$},\label{eq:NS-genAD}
\end{align}
where $\mathbf{A}_i = \frac{\partial \mathbf{F}_i   }{\partial \mathbf{U}}$ represent the Euler Jacobians, and $\mathbf{K}^h_{ij} \frac{\partial \mathbf{U}     }{\partial x_j} = \mathbf{E}_i $ the diffusivity matrices. $\mathbf{A}_i$ and $\mathbf{K}^h_{ij}$ are defined according to the set of solution variables (conservation variables in this case). For a detailed discussion on the various variable sets and corresponding $\mathbf{A}_i$ and $\mathbf{K}^h_{ij}$ matrices see, e.g., \cite{hauke1998}.
\section{Stabilized Space-Time Finite Element Formulation}
\label{sec:numerical-formulation}
Following \cite{aliabadi1995}, the deformable spatial domain/stabilized space-time (DSD/SST) finite element formulation is derived for the quasi-linear form in equation \eqref{eq:NS-genAD}.

\begin{figure}[t]
\sidecaption[t]
\includegraphics[width=0.5\textwidth]{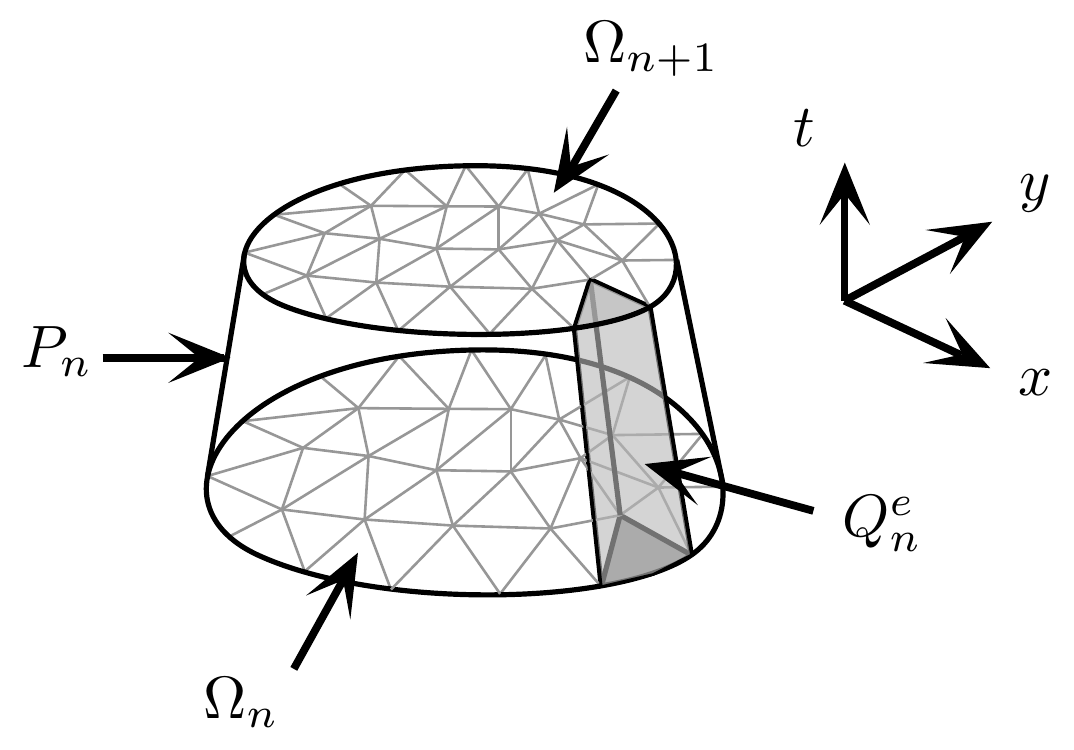}
\caption{Space-time slab with space-time element $Q_n^e$. The space-time slab is enclosed by spatial domains $\Omega_n$ and $\Omega_{n+1}$ together with $P_n$.}
\label{fig:time-slab} 
\end{figure}

In order to construct the finite element function spaces for the DSD/SST formulation, the time interval $(0,T)$ is decomposed into subintervals $I_n = (t_n, t_{n+1})$, where $t_n$ and $t_{n+1}$ are part of the ordered series: $0 = t_0 < t_1 < \cdots < t_N = T$. Additionally, we define $\Omega_n = \Omega_{t_n}$, and $\Gamma_n = \Gamma_{t_n}$. $Q_n$ now represents a so-called \textit{space-time slab}, which is the domain enclosed by $\Omega_n$, $\Omega_{n+1}$ and $P_n$. Here, $P_n$ is the surface described by the boundary $\Gamma_t$ along interval $I_n$ (see Figure \ref{fig:time-slab}). 

Similar to $\Gamma_t$ given in Section \ref{sec:governing-equations}, $P_n$ can be decomposed into $(P_n)_{g_d}$ and $(P_n)_{h_d}$. The discrete finite element space-time function spaces for the trial and weighting functions $\mathbf{U}$ and $\mathbf{W}$ are given by:
\begin{align}
(\mathcal{S}_U^h)_n &= \left\{ \mathbf{U}^h | \mathbf{U}^h \in \left[ H^{1h}(Q_n)  \right]^{n_{dof}},\: \mathbf{U}^h \cdot \mathbf{e}_d = g_d^h \:\: \mbox{on $(P_n)_{g,d}\:$, $ d= 1...n_{dof}$} \right\}, \\
(\mathcal{V}_W^h)_n &= \left\{ \mathbf{W}^h | \mathbf{W}^h \in \left[ H^{1h}(Q_n)  \right]^{n_{dof}},\: \mathbf{W}^h \cdot \mathbf{e}_d = 0 \:\:  \mbox{on $(P_n)_{h,d}\:$, $ d= 1...n_{dof}$} \right\}.
\end{align}
Using the Streamline-Upwind Petrov-Galerkin (SUPG) formulation, the weak form of Equation \eqref{eq:NS-genAD} states: given $(\mathbf{U}^h)_n^{-}$, find $\mathbf{U}^h \in (\mathcal{S}_U^h)_n$, such that $\forall \: \mathbf{W}^h \in (\mathcal{V}_W^h)_n$: 
\begin{align}
\int_{Q_n} \mathbf{W}^h &\cdot \left( \frac{\partial \mathbf{U}^h      }{\partial t}  + \mathbf{A}_i^h \cdot \frac{\partial \mathbf{U}^h  }{\partial x_i}  \right) dQ +
\int_{Q_n} \left( \frac{\partial \mathbf{W}^h  }{\partial x_i} \right) \cdot \left( \mathbf{K}^h_{ij} \frac{\partial \mathbf{U}^h  }{\partial x_j} \right) dQ + \nonumber \\
&\sum_{e=1}^{(n_{el})_n} \int_{Q^e_n}\boldsymbol{\tau}_{mom} 
\left[ (\mathbf{A}_k^h)^T  \frac{\partial \mathbf{W}^h  }{\partial x_k} \right] \cdot  
\left[ \frac{\partial \mathbf{U}^h      }{\partial t} + \mathbf{A}_i^h  \frac{\partial \mathbf{U}^h  }{\partial x_i} - \frac{\partial                  }{\partial x_i} \left( \mathbf{K}^h_{ij} \frac{\partial \mathbf{U}^h  }{\partial x_j} \right) \right]dQ +  \nonumber \\
&\sum_{e=1}^{(n_{el})_n}\int_{Q^e_n} \boldsymbol{\tau}_{DC} \left( \frac{\partial \mathbf{W}^h  }{\partial x_i} \right) \cdot \left( \frac{\partial \mathbf{U}^h  }{\partial x_i} \right) dQ + \nonumber\\
&\int_{\Omega_n} ( \mathbf{W}^h )_n^+ \cdot \left( (\mathbf{U}^h)_n^+ - (\mathbf{U}^h)_n^-   \right) d\Omega  = \int_{(P_n)_h} \mathbf{W} \cdot \mathbf{h}^h dP. \label{eq:weakform}
\end{align}
Here, the first two integrals on the left-hand side and the integral on the right-hand side represent the standard Galerkin form. The third and fourth left-hand side integrals represent the SUPG stabilisation and shock-capturing terms respectively. Continuity over the time-slab interface $\Omega_n$ is weakly imposed by the jump term, i.e., the fifth left-hand side integral in Equation \eqref{eq:weakform}. Here, the $\pm$ subscripts refer to the upper and lower time-slab solutions at time $n$. 
Equation \eqref{eq:weakform} is solved sequentially for all space-time slabs $Q_1,Q_2,\dotsc, Q_{N-1}$ with initial condition $(\mathbf{U}^h)_0^{-} = \mathbf{U}_0$. 

For $\boldsymbol{\tau}_{mom}$ in the SUPG stabilization term in Equation \eqref{eq:weakform}, the formulation proposed in \cite{hughes1986e} was used.
The shock-capturing parameter $\boldsymbol{\tau}_{DC}$ used in this work, is similar to that given in \cite{kirk2009}, which is a modification of the original definition presented in \cite{hughes1986e}. For brevity, a detailed discussion on the stabilization and shock-capturing formulations is omitted. For an extensive discussion on DSD/SST finite elements for compressible flow problems including SUPG-stabilization and shock-capturing, please refer to \cite{hughes2010}.
\section{NURBS-Enhanced Finite Elements}
\label{sec:nurbs-enhanced-finite-elements}
\begin{figure}[t]
\sidecaption[t]
\includegraphics[width=0.45\textwidth]{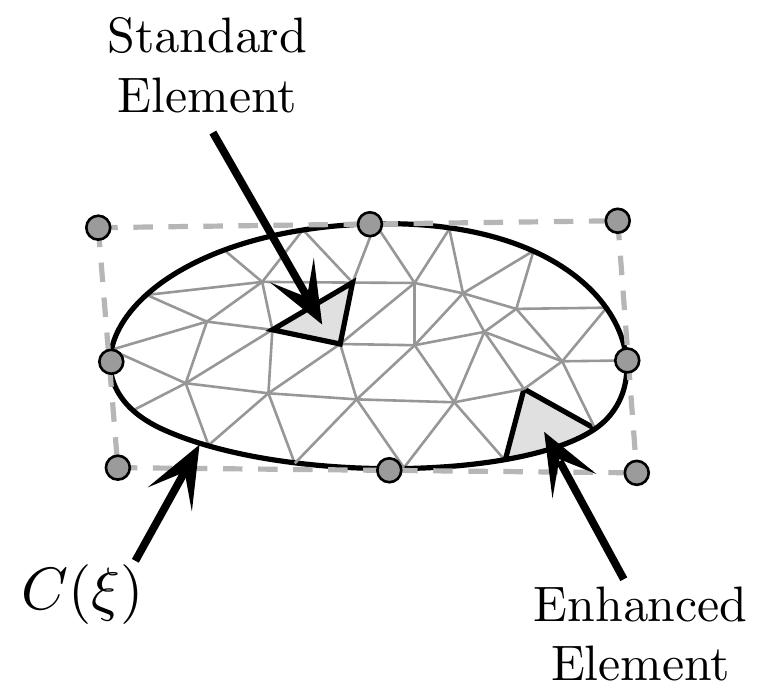}
\caption{The omputational domain defined by a NURBS curve $\mathbf{C}(\xi)$ expressed by means of parametric coordinate $\xi$ and a control polygon. The NURBS-enhanced elements are located along the NURBS boundary. The standard elements are located in the remaining interior part of the domain.}
\label{fig:nurbs-domain} 
\end{figure}

In this section, the NURBS-enhanced finite element method as proposed in \cite{hosters2018} will be presented. The key idea of this method, is to use a NURBS definition of the computational domain to \textit{enhance} the finite elements along the domain boundary. On all remaining elements in the interior of the domain a standard finite element formulation is used (cf. Figure \ref{fig:nurbs-domain}). 

Before discussing how the boundary elements make use of the NURBS boundary, let us first define a NURBS-curve of degree $p$. Such a curve is composed of piecewise rational basis functions $\mathbf{R}_i^p(\xi)$, and control points $\mathbf{B}_i$. The curve is then expressed by means of parametric coordinate $\xi \in (0,1)$ as follows:
\begin{equation}
    \mathbf{C}(\xi) = \sum_{i=1}^n \mathbf{R}_i^p(\xi) \mathbf{B}_i,
\label{eq:NURBS}
\end{equation}
where $n$ denotes the total number of control points. 

The elements that touch the NURBS domain boundary make use of a non-linear mapping between a reference element and the element in physical coordinates. This mapping, was proposed in \cite{hosters2018} as Triangle-Rectangle-Triangle (TRT) mapping. The mapping $\boldsymbol{\Phi}(s,r)$ is given by: %
\begin{equation}
\boldsymbol{\Phi}(s,r) = (1-s-r)x_3 + (s+r)\mathbf{C}\left(\frac{s\: \xi_1 +r\: \xi_2 }{s+r}\right).
\label{eq:TRT-mapping}
\end{equation}
\begin{figure}[t]
\sidecaption[t]
\includegraphics[width=0.5\textwidth]{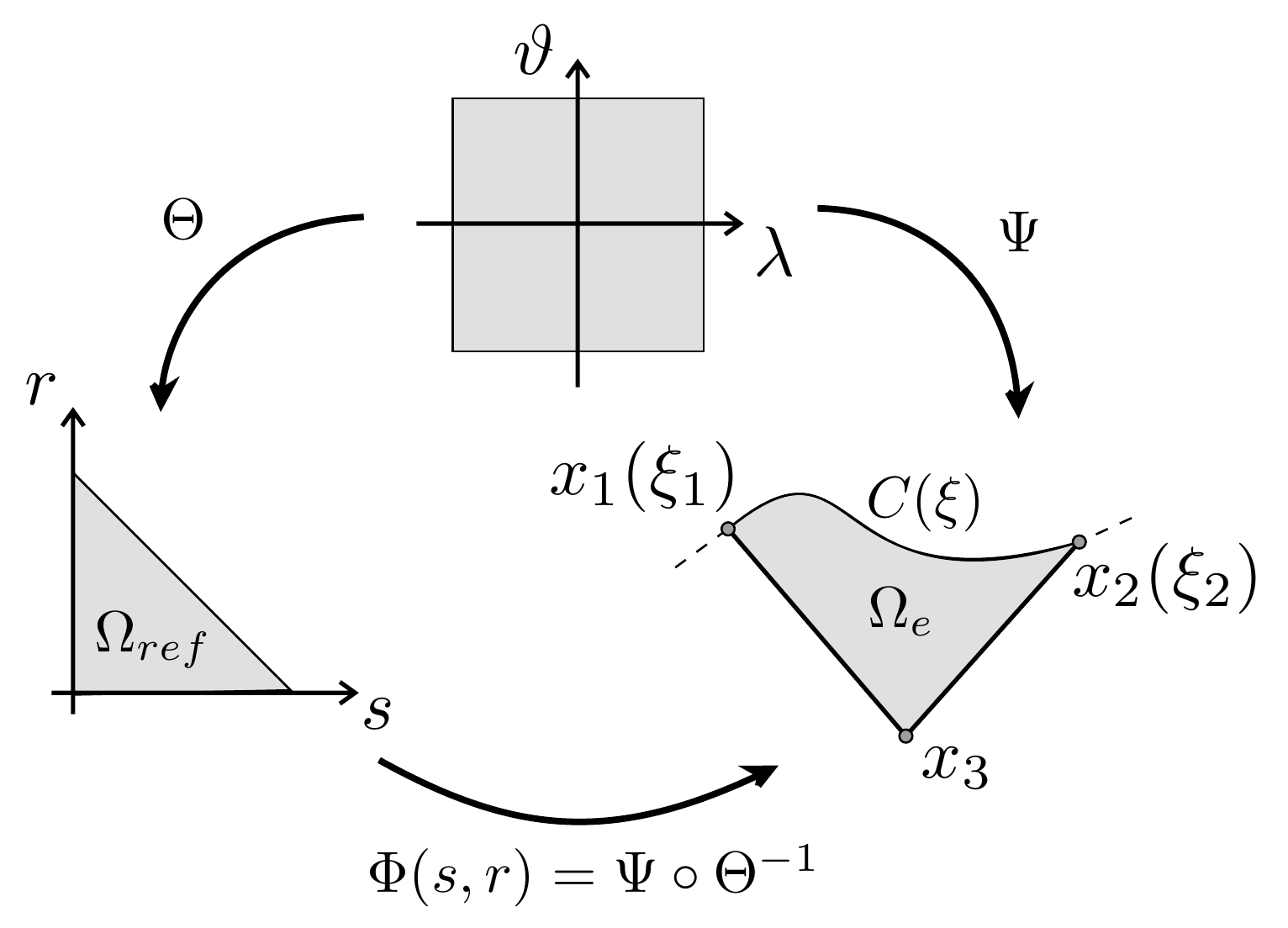}
\caption{Triangle-Rectangle-Triangle (TRT) mapping, as used in the NEFEM formulation.}
\label{fig:trt-mapping} 
\end{figure}
Here, $s$ and $r$ are the parametric coordinates of the triangular reference element, $x_3$ denotes the physical coordinate of the interior node, and $\xi_1$ and $\xi_2$ are the parametric coordinates of the NURBS curve at which the element boundary nodes are located. A graphical representation of this mapping is shown in Figure \ref{fig:trt-mapping}.

\begin{figure}[t]
\sidecaption[t]
\includegraphics[width=0.64\textwidth]{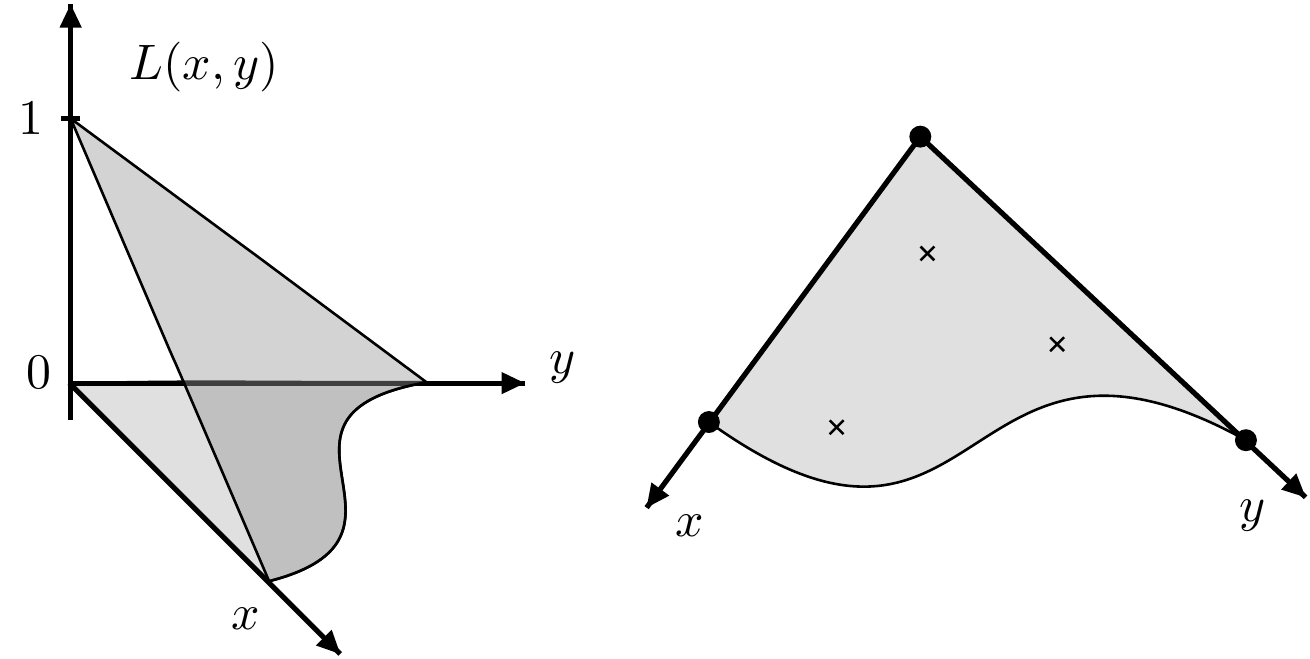}
\caption{\textbf{Left:} interior shape function, $L(x,y)$, on an NEFEM element using TRT mapping. \textbf{Right:} distribution of integration points on an NEFEM element.}
\label{fig:nefem-shapefunction} 
\end{figure}

By using the TRT mapping, the NURBS definition can be incorporated into the numerical analysis. As a result, the distribution of the integration points is determined from the exact geometry and not the erroneous discretized geometry (cf. Figure \ref{fig:nefem-shapefunction}). 

Furthermore, the shape functions corresponding to the interior nodes of the boundary elements remain zero along the NURBS (cf. Figure \ref{fig:nefem-shapefunction}). This has the advantage that there is no contribution of the interior nodes when considering Dirichlet boundaries or boundary integrals. Especially for interface-coupled problems this can be important, where Dirichlet boundaries and boundary integrals are used to compute the coupling conditions \cite{hosters2018}.   
\section{Numerical Examples}
\label{sec:numerical-examples}
To demonstrate the performance of the NEFEM in comparison to standard finite elements (SFEM), two test cases are considered next: 1) 2D supersonic viscous flow around a cylinder; 2) 2D transonic inviscid flow around a NACA0012 airfoil. 
\subsection{Cylinder Flow}
\label{sec:cylinder}
\begin{figure}[t]
\sidecaption[t]
\includegraphics[width=0.55\textwidth]{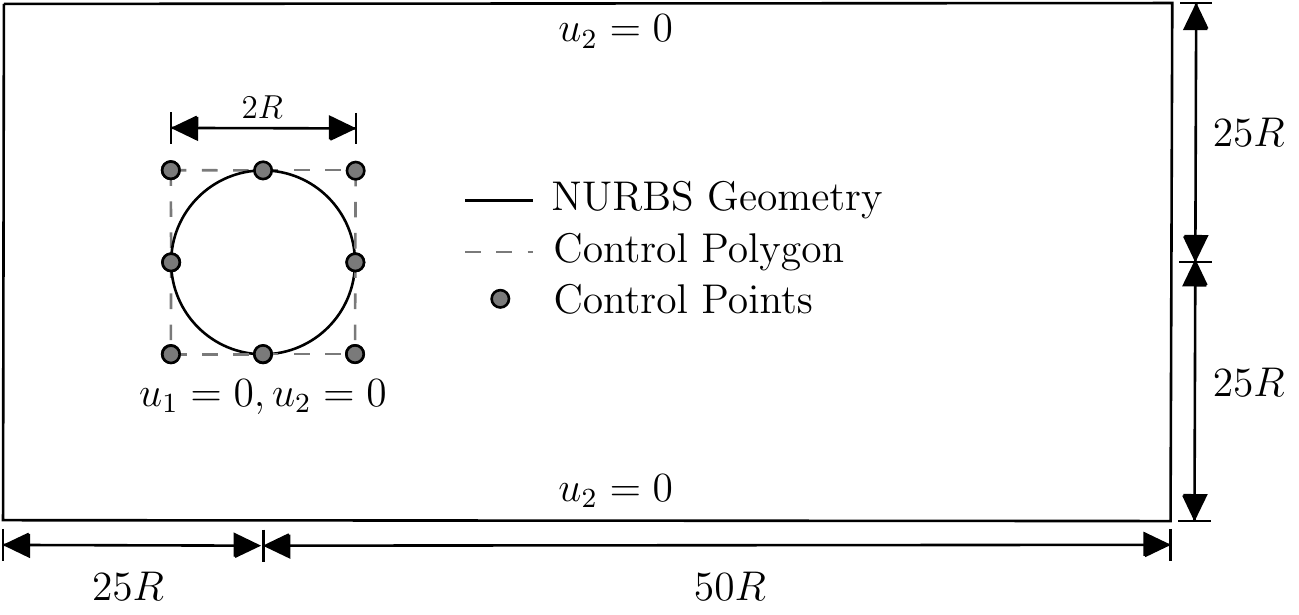}
\caption{Computational domain and boundary conditions for a supersonic flow around a 2D cylinder.}
\label{fig:cylinder-setup} 
\end{figure}

The supersonic flow around a 2D cylinder is computed using the NEFEM and the SFEM. The flow conditions and the computational domain are shown in Table \ref{tab:cylinder-setup} and Figure \ref{fig:cylinder-setup} respectively. For the NEFEM computations, the cylinder is represented by a second order NURBS-curve (cf. Figure \ref{fig:cylinder-setup}). 
\begin{table}
\caption{Flow conditions for the supersonic flow around the cylinder.}
\label{tab:cylinder-setup}
%
%
\begin{tabular}{p{1.4cm}p{1.2cm}}
\hline\noalign{\smallskip}
\multicolumn{2}{c}{Flow conditions} \\ 
\noalign{\smallskip}\svhline\noalign{\smallskip}
Mach   & 1.7             \\
Re     & $2.0\times10^5$      \\
$\rho_{in}$ & 1.0             \\
$u_{in}$    & 1.0             \\
$v_{in}$    & 0.0             \\
$e_{in}$    & 1.1179\\
\noalign{\smallskip}\hline\noalign{\smallskip}
\end{tabular}
\end{table}

\begin{figure}[t]
\sidecaption
\includegraphics[width=0.5\textwidth]{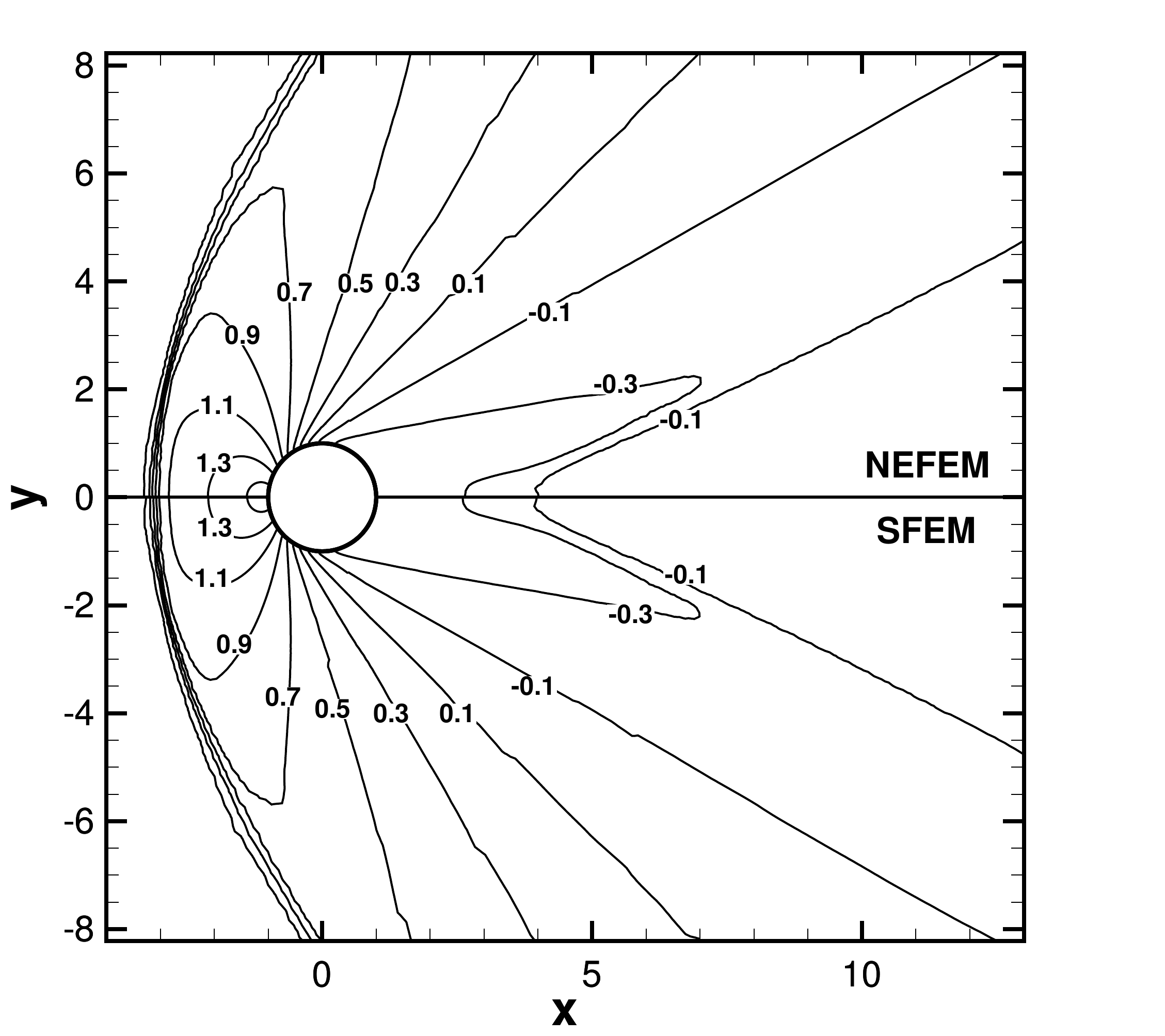}
\includegraphics[width=0.49\textwidth]{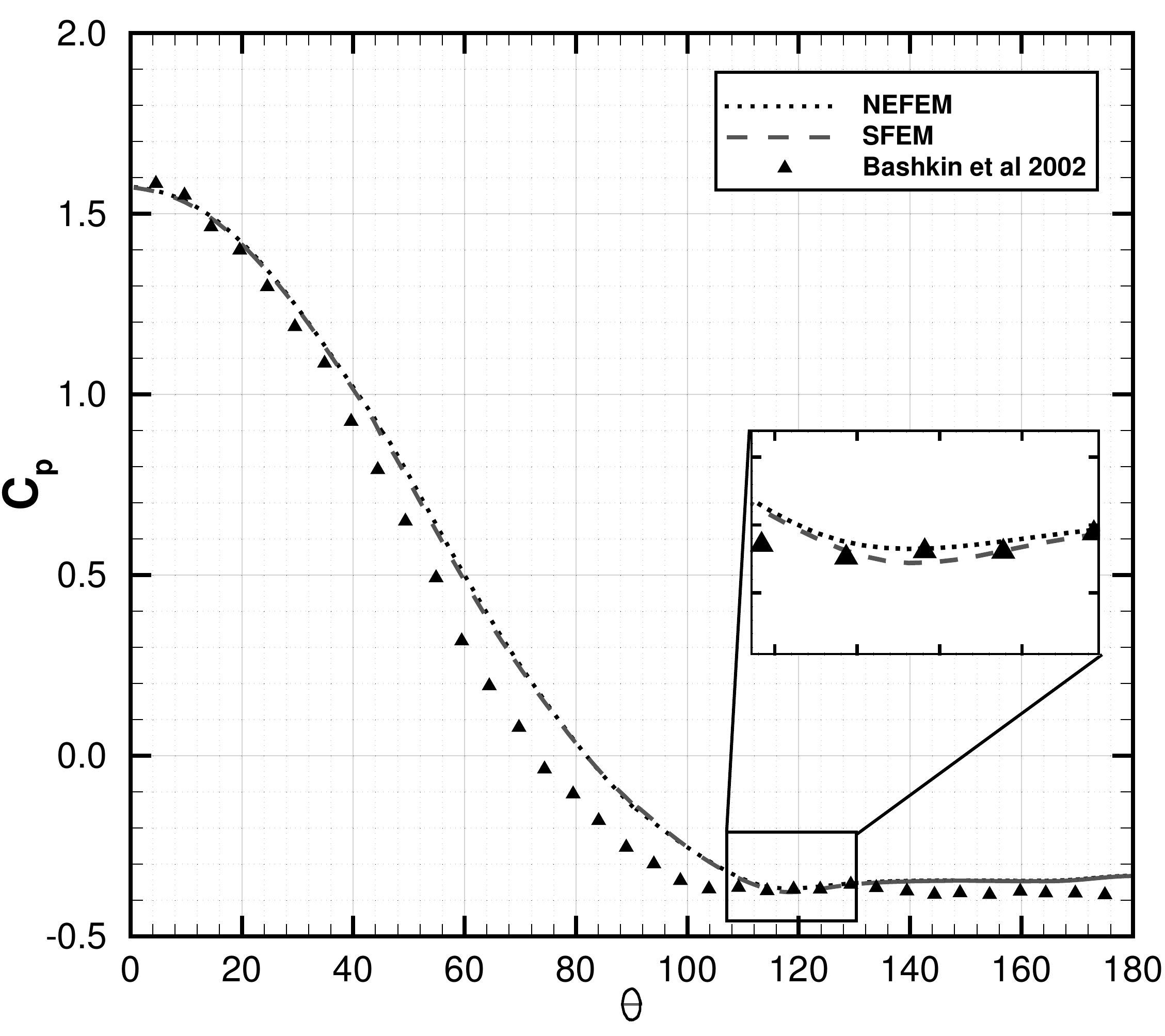}
\caption{\textbf{Left:} contour lines of the pressure coefficient for NEFEM and SFEM. \textbf{Right:} pressure coefficient along the cylinder wall with angular coordinate $\theta$.}
\label{fig:cp-flow} 
\end{figure}

The flow solution for both methods, presented by means of the pressure coefficient, $C_p$, is given in Figure \ref{fig:cp-flow}. In this figure, it can be seen that the NEFEM and the SFEM result in similar flow solutions. However, small differences can be observed when looking at the pressure coefficients along the cylinder wall. These differences could be a result of the improved geometry representation within the NEFEM. Overall, the pressure coefficient along the cylinder wall is in good agreement with the reference solution \cite{bashkin2003}, as shown in Figure \ref{fig:cp-flow}.

To demonstrate the performance of the NEFEM compared to the SFEM, a grid refinement study is performed in which the drag coefficient, $C_D$, is compared. The grids used in the study are presented in Table \ref{tab:cylinder-grids}. 

\begin{figure}[t]
\sidecaption[t]
\includegraphics[width=0.54\textwidth]{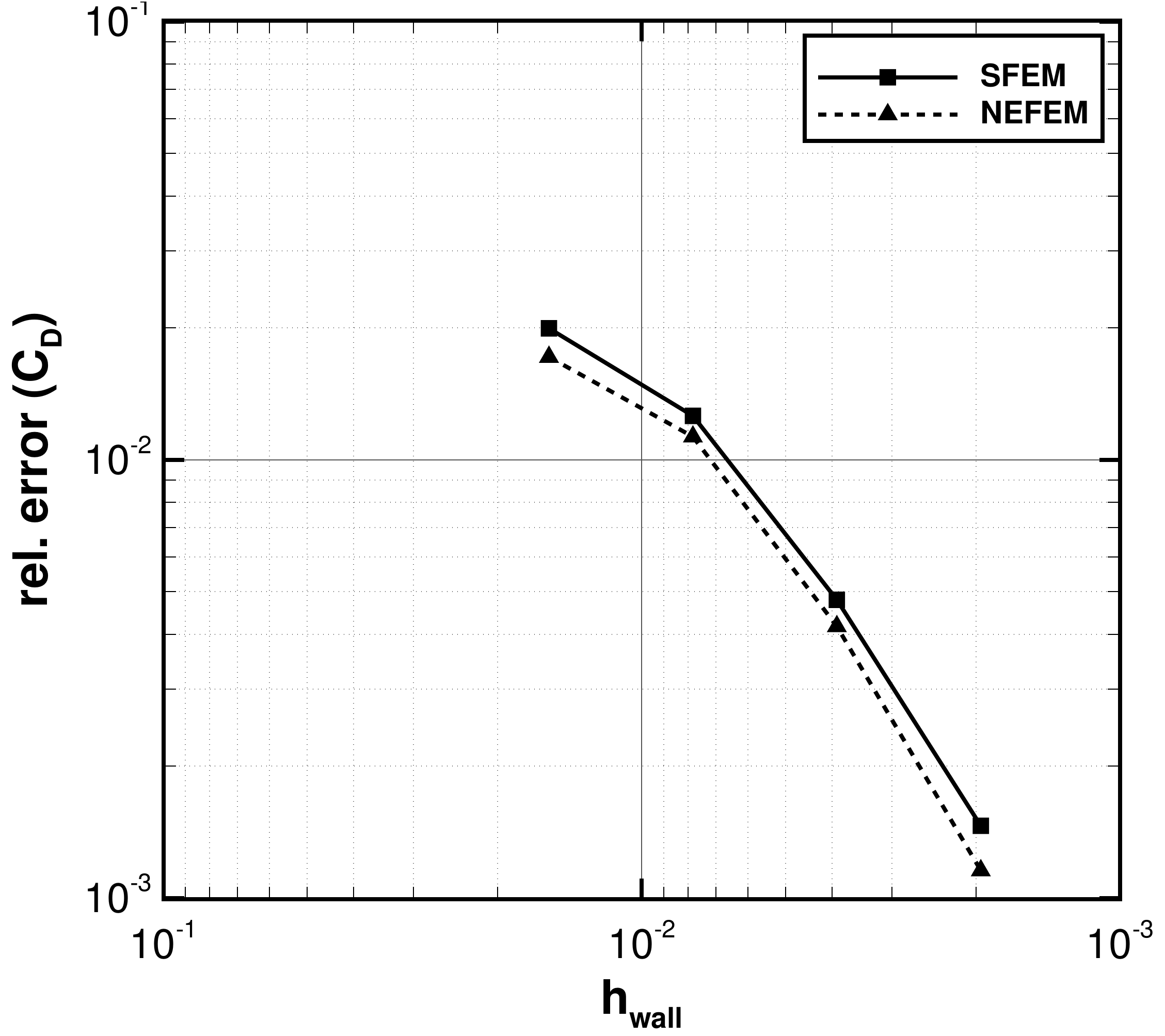}
\caption{Grid convergence for supersonic flow around a cylinder. Results are relative to that of grid 4 in Table \ref{tab:cylinder-grids}.}
\label{fig:grid-study} 
\end{figure}

The relative error in $C_D$ in Figure \ref{fig:grid-study} shows a similar convergence rate for both methods. It can be seen, however, that the NEFEM has a reduced error for all grids. 

\begin{table}
\caption{Grids used for the grid refinement study.}
\label{tab:cylinder-grids}
%
%
\begin{tabular}{p{1.0cm}p{1.4cm}p{.8cm}}
\hline\noalign{\smallskip}
Grid \# & $n_{en}$ & $n_{en_{wall}}$ \\ 
\noalign{\smallskip}\svhline\noalign{\smallskip}
0       & 6.72K    &  64             \\
1       & 26.88K   &  128             \\
2       & 107.52K  &  256            \\
3       & 430.08K  &  512            \\ 
4       & 1.72M    &  1028            \\
\noalign{\smallskip}\hline\noalign{\smallskip}
\end{tabular}
\end{table}
%
%

\subsection{NACA0012 Airfoil}
\label{sec:naca0012-airfoil}

\begin{figure}[t]
\sidecaption[t]
\includegraphics[width=0.3\textwidth]{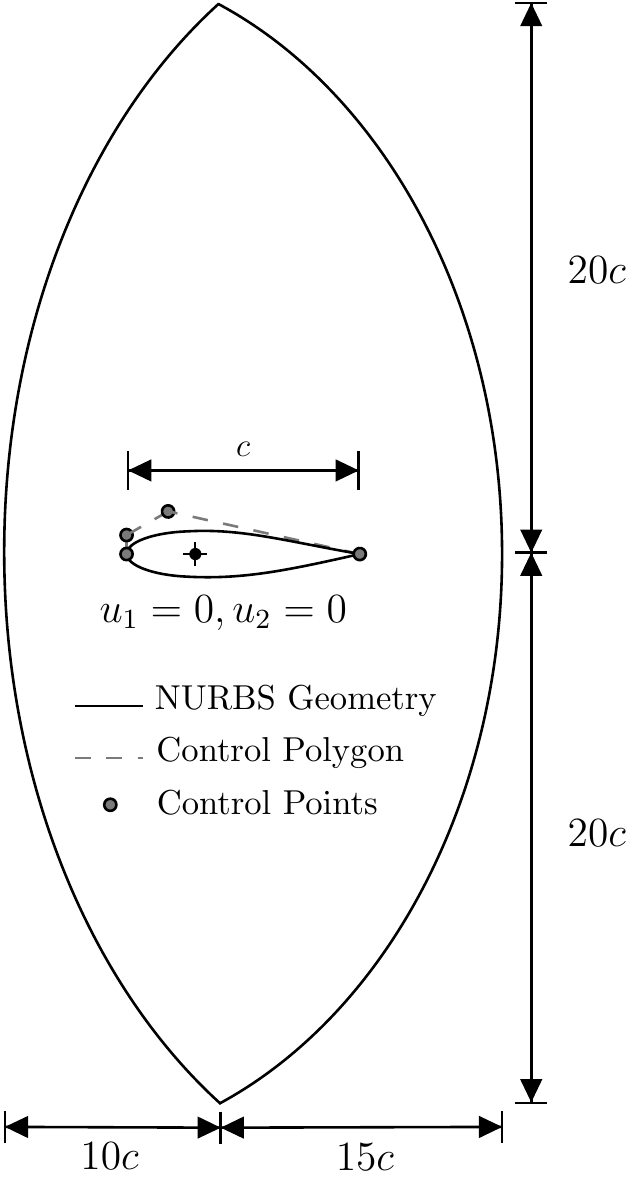}
\caption{Computational domain and boundary conditions for a supersonic flow around a NACA0012 airfoil.}
\label{fig:naca-setup} 
\end{figure}

The transonic inviscid flow around a 2D NACA0012 airfoil is computed using the NEFEM and the SFEM. The flow conditions and computational domain are shown in Table \ref{tab:naca0012-setup} and Figure \ref{fig:naca-setup} respectively. The airfoil is represented by a fourth-order NURBS curve and is positioned with a zero degree angle of attack (cf. Figure \ref{fig:naca-setup} ). 
\begin{table}
\caption{Flow conditions for the transonic flow around the NACA0012 airfoil.}
\label{tab:naca0012-setup}
%
%
\begin{tabular}{p{1.4cm}p{1.2cm}}
\hline\noalign{\smallskip}
\multicolumn{2}{c}{Flow conditions} \\ 
\noalign{\smallskip}\svhline\noalign{\smallskip}
Mach   & 0.8             \\
$\rho$ & 1.0             \\
$u$    & 1.0             \\
$v$    & 0.0             \\
$e$    & 3.29 \\
\noalign{\smallskip}\hline\noalign{\smallskip}
\end{tabular}
\end{table}

\begin{figure}[t]
\sidecaption
\includegraphics[width=0.5\textwidth]{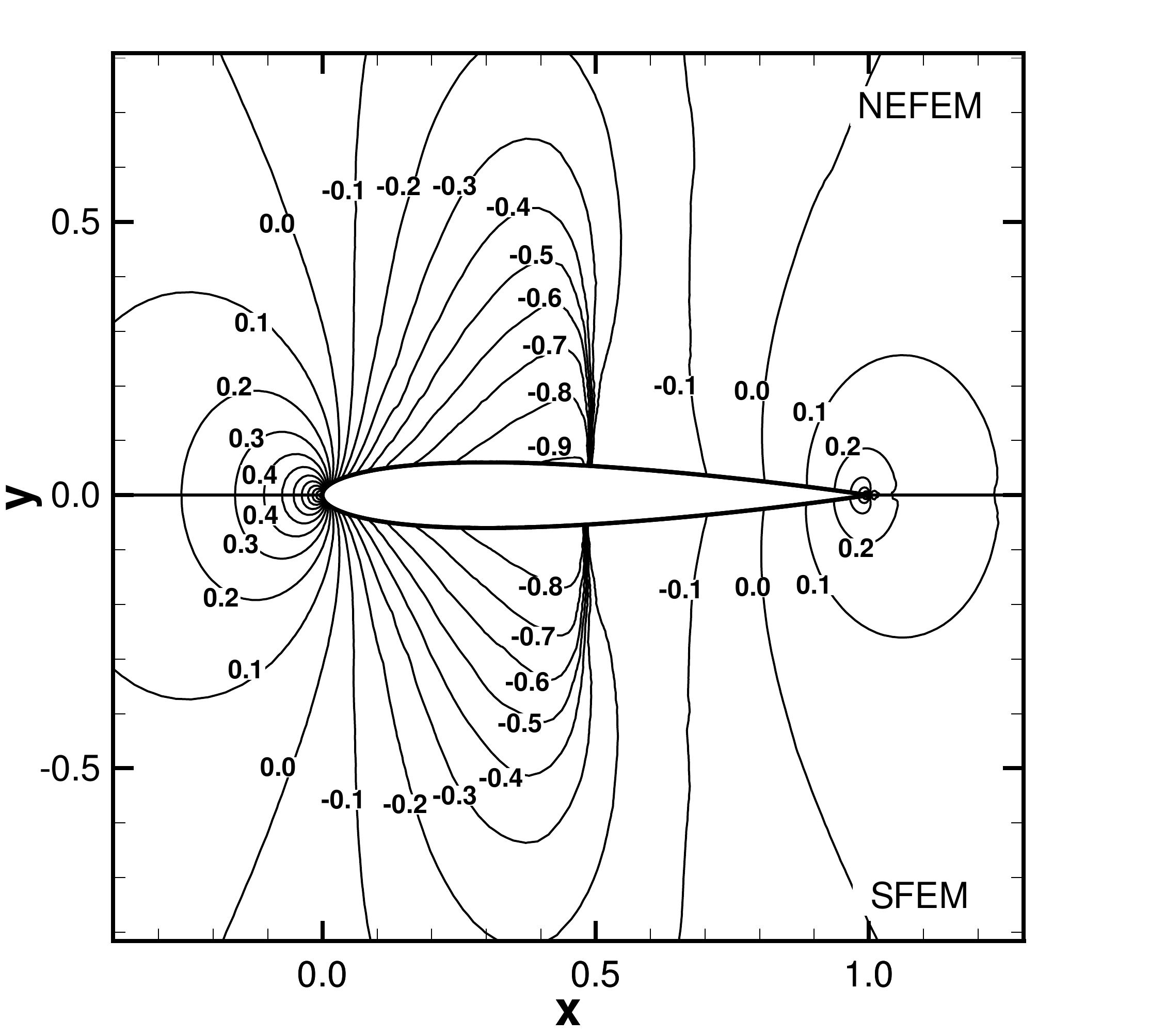}
\includegraphics[width=0.49\textwidth]{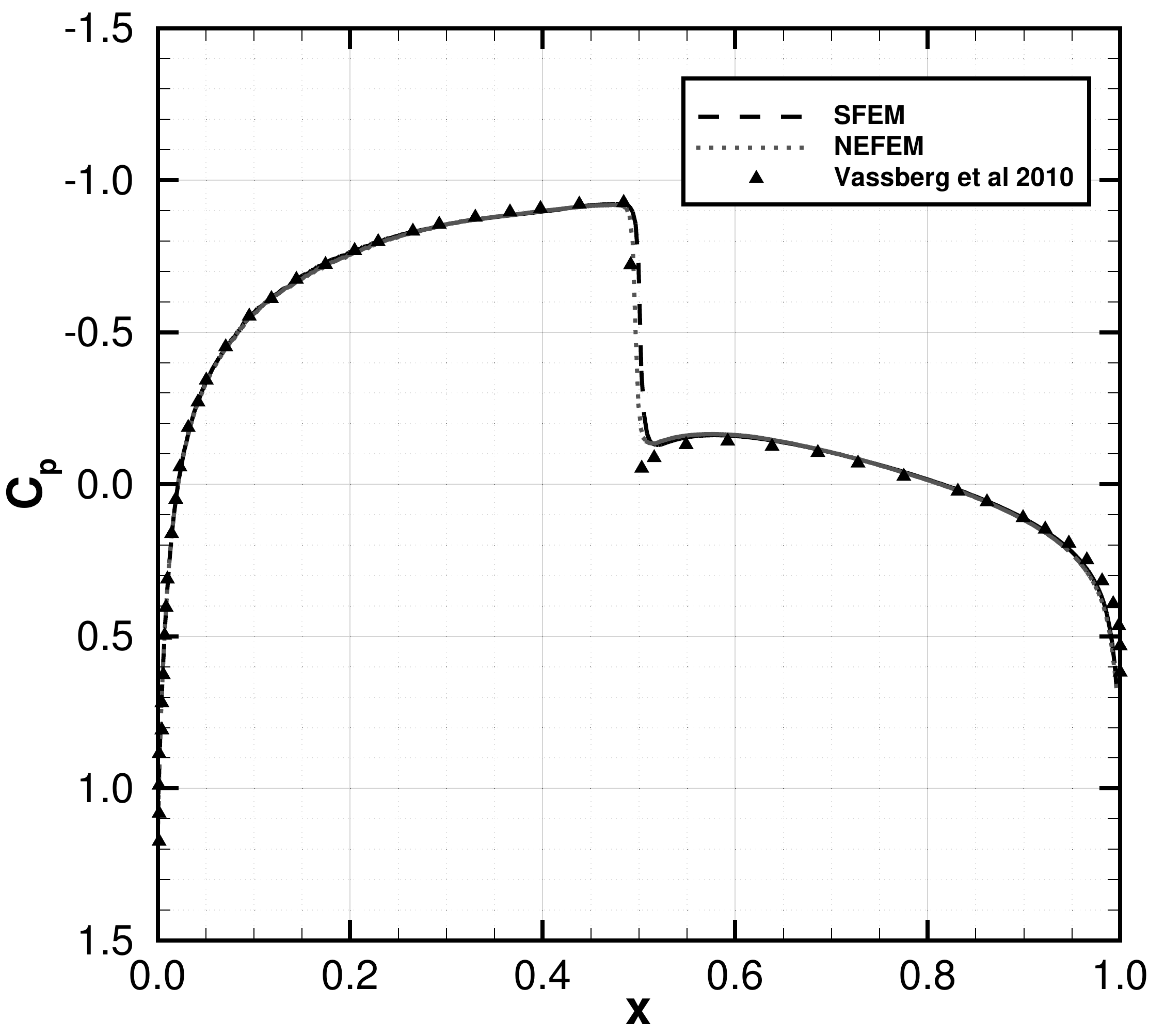}
\caption{\textbf{Left:} contour lines of the pressure coefficient $C_p$ for NEFEM and SFEM. \textbf{Right:} pressure coefficient along the airfoil wall.}
\label{fig:naca-cp-flow} 
\end{figure}

The flow solutions obtained with the NEFEM and the SFEM, are presented by means of the pressure coefficient $C_p$ in Figure \ref{fig:naca-cp-flow}. Again, it can be observed that both methods result in similar flow solutions. Additionally the results are in close agreement with the reference solution provided by \cite{vassberg2010}.

In Figure \ref{fig:naca-cp-flow}, the jump in $C_p$ along the airfoil wall is a result of the transition from supersonic to subsonic flow donditions. It can be seen that the location of this jump slightly differs between the NEFEM and the SFEM results. As for the cylinder test case, the differences between the NEFEM and the SFEM could potentially be a result of the improved geometry representation within the NEFEM.

\section{Concluding Remarks}
\label{sec:concluding-remarks}

This paper discussed the NEFEM in the context of compressible flow problems. For this method, the DSD/SST formulation was used together with the TRT mapping. Using this mapping, the NURBS definition of the domain boundary was incorporated within the shape functions and numerical integration for the elements touching these boundaries. 

The NEFEM was then tested and compared against the SFEM. Two test cases involving 2D viscous and inviscid flows were considered. Overall the methods showed good agreement with the reference solutions. Small differences between the NEFEM and SFEM solutions along solid walls in the flow domain were observed. These could be attributed to the improved geometry representation accounted for in the NEFEM. This observations could be a suitable starting point for future research on the benefits of the NEFEM over the conventional methods.

\end{document}